\documentclass[a4paper,12pt]{article}
\usepackage{latexsym}
\usepackage{amsfonts}
\usepackage[all]{xy}
\usepackage{amssymb, mathrsfs, amsfonts, amsmath}
\usepackage{amsbsy}
\usepackage{epsfig}
\newtheorem{theorem}{Theorem}[section]

\newtheorem{definition}{Definition}[section]
\numberwithin{equation}{section}

\begin{document}

\title{Non-standard eigenvalue problems for
 perturbed $p$-Laplacians}

\author{F. G\"{u}ng\"{o}r\thanks{Department of Mathematics, Faculty of Arts and Sciences,
 Do\u{g}u\c{s} University, 34722 Istanbul, Turkey, e-mail: fgungor@dogus.edu.tr and e-mail: mhasansoy@dogus.edu.tr},\and  M. Hasanov\footnotemark[1]}

\maketitle

\begin{abstract}
This paper is devoted to multi-parameter eigenvalue problems for
perturbed $p$-Laplacians, modelling travelling waves  for a class
of non-linear evolution PDE. Dispersion relations between the
eigen-para\-me\-ters, the existence of eigenvectors and positive
eigenvectors, variational principles for eigenvalues of perturbed
$p$-Laplacians and constructing  analytical solutions are the main
subject of this paper. Besides the $p$-Laplacian-like eigenvalue
problems we also deal with new and non-standard eigenvalue
problems, which can not be solved by the methods used in nonlinear
eigenvalue problems for $p$-Laplacians and similar operators. We
do both: extend and use classical variational and analytical
techniques to solve standard eigenvalue problems and suggest new
variational and analytical methods to solve the non-standard
eigenvalue problems we encounter in the search for travelling
waves.
\end{abstract}

\emph{{Keywords:}}\,\,Travelling waves, perturbed $p$-Laplacians,
eigenvalues, eigenfunctions,  variational principles, critical
points.

\emph{{AMS subject Classifications:} } Primary 49R50, 47A75,
35p15; Secondary 34L15, 35D05, 35J60

\section{Introduction}

In this paper we study a class of non-standard eigenvalue
problems, which naturally arise when we search travelling waves
for  evolution $p$-Laplacian equations. The point of departure for
the problems studied in this paper is the evolution
$(p,q)$-Laplacian equation in the following form:
$$
iv_t- \mathsf{div}(|\nabla v|^{p-2}\nabla v)= \lambda |v|^{q-2}v, \\[.3cm]
v\bigl |_{\partial Q}=0,
$$
 where $v:=v(t,x,y),\,\,t>0$, $\lambda $
is a parameter and $Q= (0,1)\times \mathbb{R}$ is an infinite
rectangle in $\mathbb{R}^2$. We can think this equation
 as a  generalized  non-linear Schr\"{o}dinger equation
because it subsumes the evolution non-linear Schr\"{o}dinger equation
in the particular case  $p=2$.

It should be noted that studying travelling waves reduces this
problem to a multi-parameter eigenvalue problem and the most
difficult eigenvalue problems arise, when $p=q$ (see  J. P. G.
Azorero, I. P. Alonso \cite{azor87}). For this reason, in this
paper we restrict ourselves to the case ${\bf p=q}$ and start with
the following problem (see \cite{bene90}):
\begin{equation}
\begin{array}{cc}
iv_t- \mathsf{div}(|\nabla v|^{p-2}\nabla v)= \lambda |v|^{p-2}v,\quad p>1 \\[.3cm]
v\bigl |_{\partial Q}=0.
\end{array}
\end{equation}
Now, we look for wave solutions to equation (1.1) in the form
$v(t,x,y)= e^{i(wt-ky)}u(x)$, where $x\in (0,1)$, $y\in
\mathbb{R}$, $u(x)$ is a real-valued function and $(w,k)\in
\mathbb{R}\times \mathbb{R}$. We have $\nabla v=
(v_x,v_y)=e^{i(wt-ky)}(u',-iku)$ and $|\nabla v|=
(k^2u^2+u'^2)^{1/2}$. For convenience we introduce the notation:
$\nabla_ku:=(ku,u')$. Then, $|\nabla v|=
|\nabla_ku|=(k^2u^2+u'^2)^{1/2}$. By using this notation and
putting $v(t,x,y)= e^{i(wt-ky)}u(x)$ into (1.1) we obtain that
$v(t,x,y)$ is a solution to equation (1.1) if and only if $u$ is
an eigenvector of the following multi-parameter eigenvalue
problem:
\begin{equation}
\begin{array}{cc}
-wu +k^2|\nabla_k u|^{p-2}u- \bigl(|\nabla_k u|^{p-2}u'\bigr)'= \lambda |u|^{p-2}u,\quad p>1 \\
u(0)=u(1)=0.
\end{array}
\end{equation}
We note that in the case of $k=0$ we obtain $|\nabla_0u|= |u'|$
and (1.2) yields
$$
-wu- (|u'|^{p-2}u')'=\lambda|u|^{p-2}u,
$$
$$
u(0)=u(1)=0,
$$
which is a two-parameter eigenvalue problem for one-dimensional
$p$-Laplacians. For this reason  we refer to  problem (1.2) as a
multi-parameter eigenvalue problem for perturbed $p$-Laplacians.
 We now give the definition
of  solution to  equation (1.2).
\begin{definition}
$0\not=u\in W_0^{1,p}(0,1)$ is a solution to  equation (1.2) if
and only if
\begin{equation}
-w\int_0^1 uv\,dx+k^2\int_0^1|\nabla_k u|^{p-2}uv\,dx+
\int_0^1|\nabla_k u|^{p-2}u'v'\,dx= \lambda\int_0^1
|u|^{p-2}uv\,dx
\end{equation}
\end{definition}
holds for all $v\in W_0^{1,p}(0,1)$, where $W_0^{1,p}(0,1)$ is the
Sobolev space (see \cite{adam02} for Sobolev spaces). Here, the
parameters $(w,k,\lambda)$ are called  eigen-parameters and the
associated non-trivial function $u\in W_0^{1,p}(0,1)$ is called
the eigenfunction. Besides the eigenfunctions, we are interested
in positive eigenfunction, too. We note that eigenvalue problems
(1.3) can be split into two groups with respect to $w$.

\textbf{I)} The case: $w\not=0$,

 \textbf{II)} The case $w=0$.

In the case \textbf{II)} problem (1.3) can be written in the
following form:
$$
 k^2\int_0^1|\nabla_k u|^{p-2}uv\,dx+ \int_0^1|\nabla_k u|^{p-2}u'v'\,dx=\lambda\int_0^1
|u|^{p-2}uv\,dx.
$$
This problem actually contains two quite different type of
eigenvalue problems:

\textbf{a)} Standard eigenvalue problems: In these problems we fix
$k$ and study $\lambda$, i.e., we are interested in the dispersion
relation $\lambda(k)$.

\textbf{b)} Non-standard eigenvalue problems: Here we search $k$
for a fixed $\lambda$.

At this point we must note that the case b) is the hardest problem
we study in this paper. Moreover, this problem is not a standard
eigenvalue problem as treated in the context of Nonlinear
Analysis. This case is the main concern of this paper and we study
it separately in the last two Sections.

First, we observe that searching for $0\not=u\in W_0^{1,p}(0,1)$
satisfying (1.3) is equivalent to finding critical points of the
functional
$$
F(u)=- \frac{w}{2}\int_0^1 u^2\,dx+\frac{1}{p}\int_0^1|\nabla_k
u|^p\,dx-\frac{\lambda}{p}\int_0^1|u|^p\,dx.
$$
In what follows, we denote  $X:=  W_0^{1,p}(0,1)$, which is a
Banach space. $F':X\rightarrow X^*$ will denote  the Fr\'{e}chet
derivative of $F$, where $X^*$ is the dual of the space $X$. It is
known that the existence of Fr\'{e}chet derivative implies the
existence of directional (Gateaux) derivative. By using the
definition of Gateaux derivative, we can obtain
\begin{eqnarray*}
&&\langle F'(u),v\rangle=\frac{d}{dt}F(u+tv)\Bigl|_{t=0}=
\frac{d}{dt}\Bigl[-\frac{w}{2}\int_0^1 (u+tv)^2\,dx+\\&&
\frac{1}{p}\int_0^1\Bigl(k^2(u+tv)^2+
(u'+tv')^2\Bigr)^\frac{p}{2}\,dx-\frac{\lambda}{p}\int_0^1|(u+tv)|^p\,dx\Bigr]\Bigl
|_{t=0}=\\&& -w\int_0^1 uv\,dx+
\frac{1}{p}\int_0^1\frac{p}{2}\Bigl(k^2u^2+
u'^2\Bigr)^\frac{p-2}{2}\Bigl(2k^2uv+2u'v'\Bigr)\,dx-\\&&
\lambda\int_0^1 |u|^{p-2}uv\,dx= -w\int_0^1
uv\,dx+k^2\int_0^1|\nabla_k u|^{p-2}uv\,dx+ \\&&\int_0^1|\nabla_k
u|^{p-2}u'v'\,dx- \lambda\int_\Omega |u|^{p-2}uv\,dx.
\end{eqnarray*}

Hence, $u\in X$ is a solution to (1.3) if and only if $u$ is a
free critical point for $F(u)$, i.e., $ \langle
F'(u),v\rangle=0,\,\,\,\hbox{ for\,\,all\,\,}v\in X $, where
$\langle F'(u),v\rangle$ denotes the value of the functional
$F'(u)$ at $v\in X$. Moreover, by Sobolev's embedding theorem  $
W_0^{1,p}(0,1)$ is compactly embedded in $C[0,1]$ and
consequently, $ W_0^{1,p}(0,1)$ is compactly embedded in
$L_q(0,1)$ for all $q\in [1,\infty)$ (see \cite{adam02}).
Therefore, the functional $F(u)$  is well defined for all $u\in
W_0^{1,p}(0,1)$.

The rest of this paper will be organized as follows.
 In Section 2 we  study the $\lambda(k)$ dependence   by using the
Ljusternik-Schnirelman critical point theory (Theorem 2.1). We
also prove a theorem (Theorem 2.3) about the existence of positive
eigenfunctions and the localization of eigen-parameters
$(w,k,\lambda)$. Section 3 is mainly devoted to constructing the
analytical solutions and the dispersion relations $k(\lambda)$
based on the analytical solutions (Theorem 3.1). Here, we modify
the existing methods and apply them to new problems. In the last
section we discuss  alternative variational methods for describing
the dispersion relations $k(\lambda)$.

\section{The structure of the eigen-parameters $(w,k,\lambda)$:\\
General results}

In this section we first study the problem
\begin{equation}
 k^2\int_0^1|\nabla_k
u|^{p-2}uv\,dx+ \int_0^1|\nabla_k u|^{p-2}u'v'\,dx=\lambda\int_0^1
|u|^{p-2}uv\,dx,
\end{equation}
and the dispersion relation $\lambda (k)$ (the case
\textbf{II)-a)}). We also present our main results about the
structure of the eigen-parameters $(w,k,\lambda)$ and associated
eigenfunctions.

 We have already noticed that for a fixed $k$ problem (2.1) is a new class of
  standard eigenvalue problem for
perturbed $p$-Laplacians. In the case of $k=0$  problem (2.1)
coincides with the typical one-dimensional eigenvalue problem for
$p$-Laplacians:
\begin{equation}
\begin{array}{cc}
-(|u'|^{p-2}u')'=\lambda|u|^{p-2}u,\,\,\,p>1 \\
u(0)=u(1)=0.
\end{array}
\end{equation}
 The eigenvalue problems for $p$-Laplacians have been studied by many
authors (see \cite{anan87}, \cite{anan96}, \cite{azor87},
\cite{le06}, \cite{lind90} and references therein). However, when
$k\not= 0$ then we deal with (2.1) which is different from the
typical one-dimensional $p$-Laplacian eigenvalue problem (2.2).
Our first observation  is given in Theorem 2.1. Since many facts
in this theorem are proved in the same way  as those in the
classical $p$-Laplacian eigenvalue problems  we present only a
sketch of the proof.
\begin{theorem} Let us fix $k\in\mathbb{R}$. Then,

a)  there exists an infinite sequence of eigenvalues
$\lambda_n(k)$ for problem (2.1), arranged as
$$
0<\lambda_1(k)<
\lambda_2(k)\leq....\leq\lambda_n(k)\leq...\,\text{and}\,\,\lambda_n(k)\rightarrow
+\infty\,\,\hbox{as}\,\,n\rightarrow\infty,
$$
where $\lambda_1(k)< \lambda_2(k)$ follows from $c)$.

 b) $ 0< \lambda_1(\Delta_p)\leq
\lambda_1(k)\,\,\hbox{and}\,\, |k|^p\leq \lambda_1(k),$

where  $\lambda_1(\Delta_p)$ denotes the first eigenvalue of
problem (2.2).

c) $\lambda_1(k)$ is simple, isolated and $\lambda_2(k)$ is the
second eigenvalue of problem (2.1)
\end{theorem}

\emph{A sketch of the Proof}. a) The proof is based on the
Ljusternik-Schnirelman critical point theory (see \cite{ze85},
Chapter 44). Let us define $G_k(u):= \frac{1}{p}\int_0^1 |\nabla_k
u|^p\,dx$ and $\Phi(u):= \frac{1}{p}\int_0^1 |u|^p\,dx$. Consider
the following eigenvalue problem:
\begin{equation}
\Phi'(u)= \mu G_k'(u),\quad u\in S_{G_k}, \quad \mu\in \mathbb{R},
\end{equation}
 where  $S_{G_k}= \{u\in X \bigl |G_k(u)=1\}$. Clearly, equation
 (2.3) is the same as equation (2.1) with $\mu=
 \frac{1}{\lambda}$. It is well known that to find a sequence $\mu_n$ of the eigenvalues
 of problem (2.3), it is sufficient to check the following basic
 conditions (see \cite{ze85}, p. 325 and p. 328):

$H1$. Let $X$ be a reflexive Banach space. $F$ and $G$ are the
even functionals such that $F,\,G\in C^1(X,\mathbb{R})$ and
$F(0)=G(0)=0$.

$H2$. $F':X\rightarrow X^*$ is strongly continuous (i.e.
$u_n\rightharpoonup u $ implies $F'(u_n)\rightarrow F'(u)$)  and
$\langle F'(u),u\rangle=0$, $u\in \overline{co}S_G$ implies $F(u)=
0$, where $\overline{co}S_G$ denotes the closure of the convex
hull of the set $S_G$.

$H3$. $G^\prime:X\rightarrow X^*$ is continuous, bounded and
satisfies the following condition:
$$
u_n\rightharpoonup u,\quad  G'(u_n)\rightharpoonup v,\quad \langle
G'(u_n),u_n\rangle\rightarrow \langle v,u\rangle
$$
implies $u_n\rightarrow u \,  $ as $n\rightarrow \infty$, where
$u_n\rightharpoonup u$ denotes the weak convergence in $X$.

$H4$. The level set $S_G$ is bounded and $u\not= 0$ implies,
$$
\langle G'(u),u\rangle>0,\quad \lim_{t\rightarrow
+\infty}G(tu)=+\infty, \quad \inf_{u\in S_G}\langle
G(u),u\rangle>0.
$$
Note that in our case, $F= \Phi$ and  $G= G_k$. It has been shown
in \cite{le06} that all conditions $H1-H4$ are satisfied for the
functionals $\Phi(u)$ and $G_0(u)$. If $k\not=0$, then one can
easily  adapt  the techniques of the proofs given in \cite{le06}
to our case, to show that all of the above-given conditions are
satisfied.

Now, we denote by $\mathcal{K}_n(k)$ the class of all compact,
symmetric subsets $K$ of $G_k$, such that ${\rm gen}\, K\geq n$
(see \cite{ze85}, Chapter 44). Thus, for a fixed $k\in
\mathbb{R}$, according to the Ljusternik-Schnirelman variational
principle (\cite{ze85}, p. 326, Theorem 44.A) there exists a
sequence of eigenvalues of problem (2.1), depending on $k$ and
arranged as:
$$
0<\lambda_1(k)\leq \lambda_2(k)\leq....\leq\lambda_n(k)\leq...
$$
 which are characterized by
$$
\frac{1}{\lambda_n(k)}=\mu_n(k)= \sup_{K\subset
\mathcal{K}_n(k)}\inf_{u\in K}\Phi(u).
$$
Moreover, for all $k\in \mathbb{R}$
$$
0<\lambda_1(k)\,\,\hbox{and}\,\, \lambda_n(k)\rightarrow
+\infty\,\,\hbox{as}\,\, n\rightarrow \infty.
$$
b) It follows from
$$
0<\lambda_1(\Delta_p)= \inf_{0\not=u\in
W_0^{1,p}(0,1)}\frac{\int_0^1|u'|^p\,dx}{\int_0^1|u|^p\,dx}\leq
\inf_{0\not=u\in
W_0^{1,p}(0,1)}\frac{\int_0^1(k^2u^2+u'^2)^{\frac{p}{2}}\,dx}{\int_0^1|u|^p\,dx}=\lambda_1(k)
$$
 and
$$
|k|^p\int_0^1|u|^p\,dx\leq \int_0^1(k^2u^2+u'^2)^{\frac{p}{2}}\,dx
$$
that $0<\lambda_1(\Delta_p)\leq \lambda_1(k)\,\,\hbox{and}\,\,
|k|^p\leq \lambda_1(k). $

c) The fact that $\lambda_1(k)$ is isolated is proved by the same
method that is  given in \cite{lind90} to prove that the first
eigenvalue of the $p$-Laplacian is isolated. Finally, to prove
that $\lambda_2(k)$ is the second eigenvalue of problem (2.1) one
may follow the approach of the paper \cite{anan96}. The simplicity
result repeats the arguments from \cite{anan87}.

\textbf{The case: I).}

In this case $w\not=0$ and our the main concern is problem (1.3)
with $u>0$:
$$
 k^2\int_0^1|\nabla_k u|^{p-2}uv\,dx+ \int_0^1|\nabla_k u|^{p-2}u'v'\,dx-\lambda\int_0^1
|u|^{p-2}uv\,dx= w\int_0^1 uv\,dx,
$$
$$
u>0\,\,\,\,\,\hbox{in}\,\,(0,1).
$$
  By using the scaling property, a solution to problem (1.3) can be obtained by a
constrained minimization problem for the functional
$$
E_{k,\lambda}(u)= \int_0^1|\nabla_k
u|^p\,dx-\lambda\int_0^1|u|^p\,dx
$$
on the Banach space $W_0^{1,p}(0,1)$, restricted to the set
$$
M=\bigl\{u\in W_0^{1,p}(0,1)\Bigl | \int_0^1u^2\,dx=1\bigr\}.
$$
The main idea is based on some regularity ideas, on  the fact that
$W_0^{1,p}(0,1)$ is compactly embedded in $L_q(0,1)$ for all $q\in
[1,\infty)$ and on the following theorem.
\begin{theorem} (see \cite{stru02},\,\,Theorem 1.2) Suppose $X$ is a reflexive
Banach space with norm $\|.\|$, and let $M\subset X$ be a weakly
closed subset of $X$. Suppose $E: M\rightarrow \mathbb{R}\cup
+\infty$ is coercive on $M$ with respect to $X$, that is

1) $E(u)\rightarrow \infty$ as $\|u\|\rightarrow\infty, \,\,u\in
M$

and it is (sequentially) weakly lower semi-continuous on $M$ with
respect to $X$, that is

2) for any $u\in M$, any sequence $(u_n)$ in $M$ such that
$u_n\rightharpoondown u$ (weakly) in $X$ there holds:
$$
E(u)\leq\lim_{n\rightarrow\infty}\inf E(u_n).
$$
Then $E$ is bounded from below on $M$ and attains its infimum in
$M$.
\end{theorem}

Now we formulate and prove our  main result in this section.
\begin{theorem}
 If either $(w,k,\lambda)\in \mathbb{R_+}\times \mathbb{R}\times
 (-\infty,\lambda_1(k))$ or $(w,k,\lambda)\in \mathbb{R_-}\times \mathbb{R}\times
 (\lambda_1(k), +\infty)$, then problem (1.3) has a positive
 solution.
\end{theorem}
Proof. Let us consider the condition
$(w,k,\lambda)\in\mathbb{R_+}\times
\mathbb{R}\times(-\infty,\lambda_1(k))$. By this condition we have
to prove that for a fixed $k\in \mathbb{R}$ and
$\lambda<\lambda_1(k)$ problem (1.3) has a positive solution for
any $w>0$. We set in Theorem 2.2: $X=W_0^{1,p}(0,1)$, $E(u):=
E_{k,\lambda}(u)$ and $ M=\bigl\{u\in W_0^{1,p}(0,1)\Bigl |
\int_0^1u^2\,dx=1\bigr\} $. Evidently, all conditions of Theorem
2.2 are satisfied. Particularly, by the Sobolev's embedding
theorem $M$ is a weakly closed set. Now, the existence of a
non-trivial solution to problem (1.3) immediately follows from
this theorem. The existence of a non-negative solution is obtained
if we replace $u$ by $|u|$. To prove that a non-negative solution
is positive we use some regularity results for solutions to (1.3).
We can do this in the three steps given below.

\emph{Step 1.} We show that a solution to (1.3) belongs to
$L^\infty(0,1)$. To show this one can use the Moser iteration
technique (see \cite{drab97} or \cite{le06}, pp. 1070-1073).
Actually, one can repeat step by step the method which was applied
to prove that the eigenfunctions for $p$-Laplacians are bounded
(see \cite{le06}).

\emph{Step 2.} Now, we prove that $u\in C^{1,\alpha}(0,1)$-H\"{o}lder
continuously differentiable function with the exponent
$0\leq\alpha\leq 1$. The proof is based on the following fact: Let
$\Omega$ be a bounded domain in $\mathbb{R}^n$ and
$f:\Omega\times\mathbb{R}\rightarrow \mathbb{R}$ be a
Carath\'{e}dory function (see \cite{le06}, p. 1074).  Then if
$g(x):= f(x,u(x))\in L^\infty(0,1)$ then a result of DiBenedetto
\cite{bene83} and Tolksdorf \cite{tolk84}  states that a weak
solution of the equation
\begin{equation}
-\Delta_pu(x)= f(x,u(x))\,\,\,\hbox{in}\,\,\Omega
\end{equation}
is a $C^{1,\alpha}(\Omega)$ function.

\emph{Step 3.} Finally, we use the following Harnack type
inequality due to Trudinger (\cite{gilb01} and \cite{le06}, p.
1075) to prove that $u>0$.

\emph{Harnack inequality:}  Let $u\in W^{1,p}(\Omega)$ be a weak
solution of (2.4) and for all $M<\infty$ and for all $(x,s)\in
\Omega\times (-M,M)$ the condition
$$
|f(x,s)|\leq b_1(x)|s|^{p-1}+b_2(x)
$$
holds, where $b_1$, $b_2$ are nonnegative functions in
$L^\infty(\Omega)$. Then if $0\leq u(x)<M$ in a cube
$K(3r):=K_{x_0}(3r)\subset \Omega$, there exists a constant $C$
such that
$$
\max_{K(r)}u(x)\leq C\min_{K(r)}u(x).
$$
In our problem it is enough to check \emph{Step 1}. Then the chain
\emph{Step 1}$\Rightarrow$ \emph{Step 2} $\Rightarrow$ \emph{Step
3} is obvious. We also note that the condition $
\max_{K(r)}u(x)\leq C\min_{K(r)}u(x) $
 means either: $u=0$ or $u>0$ in $\Omega$. Since
$\|u\|_{L_2(0,1)}=1$  we obtain that $u>0$ in $\Omega$.

The rest of the paper is devoted to case \textbf{II-b)}, i.e., we
seek $k$ for a fixed $\lambda$ in problem (2.1).

\section{Analytical solutions of two-parameter eigenvalue problems and dispersion relations}

In this section we study the analytical solutions of the following
problem:

\begin{equation}
\begin{array}{cc}
k^2|\nabla_k u|^{p-2}u- \bigl(|\nabla_k u|^{p-2}u'\bigr)'= \lambda |u|^{p-2}u,\quad p>1 \\
u(0)=u(1)=0.
\end{array}
\end{equation}

It turns out that in some cases it is possible to find an
analytical solution to problem (3.1). The construction of the
analytical solutions also allows us to get some dispersion
relations between $k$ and $\lambda$. In the next Section, we will
separately discuss  the methods of describing the dispersion
relations $k(\lambda)$ in the case when we can not find analytical
solutions.

First we note that in the case $k=0$ we deal with the following
classical one-dimensional $p$-Laplacian eigenvalue problem which
was fully studied by P. Dr\'{a}bek \cite{drab80} (see also a paper
of M. Del Pino, M. Elgueta, R. Manasevich \cite{pino89}):

\begin{equation}
\begin{array}{cc}
-(|u'|^{p-2}u')'=\lambda|u|^{p-2}u,\quad p>1 \\
u(0)=u(1)=0.
\end{array}
\end{equation}
All eigenfunctions and eigenvalues are given by $u_n(x)=
c\lambda_n^{-\frac{1}{p}}\sin_p(\lambda_n^{\frac{1}{p}}x)$ and\\
$\lambda_n= \bigl(n\pi_p)^p$, respectively. Here, $\pi_p:=
2\int_0^{(p-1)/p}\frac{ds}{\bigl(1-\frac{s^p}{p-1}\bigr)^{1/p}}$
and $\sin_p(x)$  is defined as an implicit function $\sin_p: [0,
\frac{\pi_p}{2}]\rightarrow [0, (p-1)^{1/p}]$ by
$$
\int_0^{\sin_p(x)}\frac{ds}{\bigl(1-\frac{s^p}{p-1}\bigr)^{1/p}}=x,
$$
then it is extended by setting: $\widetilde{\sin}_p(x):=
\sin_p(\pi_p-x),\,\,x\in [\frac{\pi_p}{2}, \pi_p]$ and
$\widetilde{\sin}_p(x):= - \widetilde{\sin}_p(-x)$ for $x\in
[-\pi_p,0]$. Finally, $\sin_p:\mathbb{R}\rightarrow \mathbb{R}$ is
defined as the $ 2\pi_p$-periodic extension of $\tilde{\sin}_p(x)$
to all of $\mathbb{R}$ (see \cite{drab80} and \cite{pino89}, and
also the recent paper \cite{taka10} for more interesting
properties of $\sin_p(x)$).

To construct analytical solutions to problem (3.1) for some
special cases, we also apply methods  similar to those applied in
the above mentioned papers. However, we have to modify some
techniques of these papers which are not applicable to our
problems. Next, we present a modified version of the methods, used
in \cite{drab80} and \cite{pino89} to construct the analytical
solutions to problem (3.2). Namely, this modified method will be
applied to solve analytically problem (3.1) in some special cases.

Let us consider equation (3.2). For the sake of simplicity we
assume that $u$ and $u'$ are positive. Thus we consider the
following problem:
\begin{equation}
\begin{array}{cc}
-((u')^{p-1})'=\lambda u^{p-1},\quad p>1 \\
u(0)=u(1)=0.
\end{array}
\end{equation}
or equivalently
$$
\begin{array}{cc}
-(p-1)(u')^{p-2}u''=\lambda u^{p-1},
 \\
u(0)=u(1)=0.
\end{array}
$$

By using the substitution $u'=v(u)$ we can reduce the order of
this equation. Indeed, $u'=v(u)$ implies $u''=vv'$. Then we have
$$
-(p-1)v^{p-1}v'=\lambda u^{p-1},
$$

$$
v^{p-1}dv=- \frac{\lambda}{p-1}u^{p-1}du.
$$
Integrating both sides we obtain
\begin{equation}
v^p=c-\frac{\lambda}{p-1}u^p \Leftrightarrow
u'^p=c-\frac{\lambda}{p-1}u^p.
\end{equation}
Now, by using the condition $u(0)=0$ we can write $
\int_0^{u(x)}\frac{dr}{\bigl(c-\frac{\lambda}{p-1}r^p\bigr)^{1/p}}=x
$ and the substitution $ s=c^\frac{1}{p}r$ yields
\begin{equation}
\int_0^{c^{-\frac{1}{p}}u(x)}\frac{ds}{\bigl(1-\frac{\lambda}{p-1}s^p\bigr)^{1/p}}=x.
\end{equation}
\textbf{Note.} At this point we have to note that, we can take
 $\lambda$ out of the integrand by using the substitution $s=s_0\lambda^{-\frac{1}{p}}$
 and then  define the function $\sin_p(x)$. The authors in \cite{drab80} and \cite{pino89}
 follow this way. However,  we will see below that,
 for problem (3.1) we also have a similar expression where we can
 not take $\lambda$ out. That is why we need a modification of
 this method.

Now we apply a technique which will also be applied to solve more
difficult  equation (3.1). Let us define
$$
F(\lambda,x):=
\int_0^{x}\frac{ds}{\bigl(1-\frac{\lambda}{p-1}s^p\bigr)^{1/p}}.
$$
Evidently,
$$F(\lambda, .): \Bigl[0,
\bigl(\frac{p-1}{\lambda}\bigr)^{\frac{1}{p}}\Bigr]\rightarrow
\Bigl[0, \frac{\pi_p(\lambda)}{2}\Bigr],\,\,\hbox{
where}\,\,\pi_p(\lambda):=
2\int_0^{(\frac{p-1}{\lambda})^{1/p}}\frac{ds}{\bigl(1-\frac{\lambda}{p-1}s^p\bigr)^{1/p}}.
$$
We have $F(\lambda,0)=0$ and $F'_x(\lambda,x)=
\frac{1}{\bigl(1-\frac{\lambda}{p-1}x^p\bigr)^{1/p}}>0,\,x\in [0,
(\frac{p-1}{\lambda})^{1/p}).$ It means that there exists the
inverse function $G(\lambda, .): [0,
\frac{\pi_p(\lambda)}{2}]\rightarrow
[0,(\frac{p-1}{\lambda})^{\frac{1}{p}}]$ defined by
$$
\int_0^{G(\lambda,
x)}\frac{ds}{\bigl(1-\frac{\lambda}{p-1}s^p\bigr)^{1/p}}=x.
$$
To extend the function $G(\lambda,x)$ we follow the same way,
which has been applied to extend  $\sin_p(x)$. Namely,
$\tilde{G}(\lambda,x):= G(\lambda, \pi_p(\lambda)-x),\,\,x\in
[\frac{\pi_p(\lambda)}{2}, \pi_p(\lambda)]$ and
$\tilde{G}(\lambda,x)= - \tilde{G}(\lambda,-x)$ for $t\in
[-\pi_p(\lambda),0]$. Finally, $G(\lambda,x):\mathbb{R}\rightarrow
\mathbb{R}$ is defined as the $ 2\pi_p(\lambda)$-periodic
extension of $\tilde{G}(\lambda,x)$ to all $\mathbb{R}$. It
follows from this construction that, $G(\lambda,x)=0
\Leftrightarrow x= n\pi_p(\lambda),\,\,n=0,\pm 1,\pm 2,...$.
Moreover, we obtain from (3.5) that $u(x) = c G(\lambda, x)$
verifies  the equation $-((u')^{p-1})'=\lambda u^{p-1}$ and the
initial condition $u(0)=0$. Finally, to get a solution of (3.3) we
use the condition $u(1)=0$. Now, $u(1)=0\Rightarrow
G(\lambda,1)=0\Rightarrow n\pi_p(\lambda)=1\Rightarrow
\pi_p(\lambda)= \frac{1}{n}.$ Clearly, by using the substitution $
s= \frac{r}{\lambda^{1/p}}$ we obtain  $\pi_p(\lambda)=
\frac{\pi_p(1)}{\lambda^{1/p}}$. So  $\pi_p(\lambda)= \frac{1}{n}$
implies $\lambda_n= (n\pi_p(1))^{1/p}$. Since this argument will
not work for problem (3.1), we modify this in the following way:
we avoid constructing an exact solution for the equation
$\pi_p(\lambda)= \frac{1}{n}$ and instead we show that the
eigenvalues of problem (3.3) consist of a sequence $\lambda_n$ and
$\lambda_n\rightarrow +\infty$ as $n\to \infty$. Indeed, we have
$\lim_{\lambda\rightarrow 0+}\pi_p(\lambda)=+\infty$,
$\lim_{\lambda\rightarrow +\infty}\pi_p(\lambda)=0$ and
$\pi_p(\lambda)$ is a decreasing function on $(0,+\infty)$. It
follows from these properties that, for each $n\in \mathbb{N}$ the
equation $\pi_p(\lambda)= \frac{1}{n}$ has a unique solution
$\lambda_n$ and $\lambda_n\rightarrow +\infty$ as $n\rightarrow
\infty$.

Now, we demonstrate  these modified techniques on the following
model problem to get analytical solutions and dispersion
relations.

\subsection{A Model Problem:  $p=4$ and $k\not=0$.}

In this Subsection we set $p=4$ in (3.1) and try to find its
analytical solution and some dispersion relations between $k$ and
$\lambda$ for the following equation:

\begin{equation}
\begin{array}{cc}
k^2(k^2u^2+u'^2)u- \bigl((k^2u^2+u'^2)u'\bigr)'=\lambda u^3 \\
u(0)=u(1)=0,
\end{array}
\end{equation}
where $u=u(x)$. In what follows we also assume $k\geq 0$ because
of the symmetry property with respect to $k$.

 First, we reduce the order of the equation by the
substitution $u'=v(u)$. Then, $u''=v'=
\frac{dv}{du}\frac{du}{dx}=\frac{dv}{du}v$ and we have
%$$
%k^2(k^2u^2+u'^2)u-\bigl[
%(k^2u^2+u'^2)'u'+(k^2u^2+u'^2)u''\bigr]=\lambda u^3,
%$$
%$$
%k^2(k^2u^2+u'^2)u-
%\bigl[2k^2uu'^2+2u'^2u''+k^2u^2u''+u'^2u''\bigr]=\lambda u^3,
%$$
%$$
%k^4u^3+k^2uv^2-2k^2uv^2-2v^3\frac{dv}{du}-k^2u^2v\frac{dv}{du}-v^3\frac{dv}{du}=\lambda
%u^3,
%$$
$$
(3v^3+k^2u^2v)\frac{dv}{du}=k^4u^3-k^2uv^2-\lambda u^3.
$$

This is a homogeneous ordinary differential equation and it may be
integrated by changing from $(u,v)\to (u,w)$
%\begin{equation}
%\Bigl(3\bigl(\frac{v}{u}\bigr)^3+k^2\bigl(\frac{v}{u}\bigr)\Bigr)\frac{dv}{du}=k^4-k^2\bigl(\frac{u}{v}\bigr)^2-\lambda.
%\end{equation}
by the standard substitution $v=uw(u)$. This then gives a
separable equation
%$$
%(3w^3+k^2w)(w+uw')=k^4-\lambda-k^2w^2,
%$$
%$$
%3w^4+k^2w^2+ (3w^3+k^2w)uw'=k^4-\lambda-k^2w^2,
%$$
$$
(3w^3+k^2w)uw'=k^4-\lambda-2k^2w^2-3w^4,
$$
%$$
%\frac{3w^3+k^2w}{k^4-\lambda-2k^2w^2-3w^4}dw= \frac{du}{u},
%$$
%
%$$
%\int\frac{3w^3+k^2w}{k^4-\lambda-2k^2w^2-3w^4}dw=
%\int\frac{du}{u},
%$$
%$$
%-\frac{1}{4}\int\frac{d(k^4-\lambda-2k^2w^2-3w^4)}{k^4-\lambda-2k^2w^2-3w^4}=\int\frac{du}{u},
%$$
%$$
%\ln|k^4-\lambda-2k^2w^2-3w^4|=-4\ln|u|+c,
%$$
%
%$$
%\ln\bigl(u^4|k^4-\lambda-2k^2w^2-3w^4|\bigr)=c
%$$
which integrates to
$$
k^4-\lambda-2k^2w^2-3w^4= \frac{c}{u^4}.
$$
Putting the inverse substitution $w=\frac{u'}{u}$ into the above
equation we obtain
$$
k^4-\lambda-2k^2\Bigl(\frac{u'}{u}\Bigr)^2-3\Bigl(\frac{u'}{u}\Bigr)^4=
\frac{c}{u^4}
$$
or
$$
3u'^4+2k^2u^2u'^2-(k^4-\lambda)u^4-c=0
$$
through the multiplication by $-u^4$. Finally, by solving this
quadratic equation in $u'^2$ for $u'$ we get
%$$
%u'^2=
%\frac{-k^2u^2+\sqrt{k^4u^4+3\bigl((k^4-\lambda)u^4+c\bigr)}}{3}
%$$
%and
$$
u'=
\Bigl[-\frac{k^2}{3}u^2+\sqrt{\frac{k^4}{9}u^4+\frac{1}{3}\bigl((k^4-\lambda)u^4+c\bigr)}
\,\,\,\,\,\Bigr]^{1/2}.
$$
We notice that if $k=0$ and $p=4$ then this equation and equation
(3.4)  coincide. Hence, integrating the last equation and using
the initial condition $u(0)=0$ gives
$$
\int_0^{u(x)}\frac{dr}{\Bigl[-\frac{k^2}{3}r^2+\sqrt{\frac{k^4}{9}r^4+\frac{k^4-\lambda}{3}r^4
+C}\,\,\,\Bigr]^{1/2}}=x.
$$
Let $r=C^{1/4}s$. Then we have

\begin{equation}
\int_0^{cu(x)}\frac{ds}{\Bigl[-\frac{k^2}{3}s^2+\sqrt{\frac{k^4}{9}s^4+
1-\frac{\lambda-k^4}{3}s^4}\,\,\,\Bigr]^{1/2}}=x.
\end{equation}
Again, in the case of $k=0$ we obtain from this equation that
$$
\int_0^{cu(x)}\frac{ds}{\Bigl[1-\frac{\lambda}{3}s^4\Bigr]^{1/4}}=x,
$$
which is the same as (3.5). Let us define
\begin{equation}
F_{\lambda,k}(x):=
\int_0^{x}\frac{ds}{\Bigl[-\frac{k^2}{3}s^2+\sqrt{\frac{k^4}{9}s^4+
1-\frac{\lambda-k^4}{3}s^4}\,\,\,\Bigr]^{1/2}}.
\end{equation}
By Theorem 2.1 we have $|k|^p\leq\lambda_1(k)$. Hence, $k^4\leq
\lambda_1(k)$. For this reason in what follows we suppose that $k^
4<\lambda$ (the case $\lambda=k^4$ is considered in Theorem 3.1).
Function $F_{\lambda,k}(x)$ is well defined if
$1-\frac{\lambda-k^4}{3}s^4>0$. Therefore, $
s<\Bigl(\frac{3}{\lambda-k^4}\Bigr)^{1/4}$ and
$$
F_{\lambda,k}: \Bigl[0,
\Bigl(\frac{3}{\lambda-k^4}\Bigr)^{1/4}\Bigr]\rightarrow \Bigl[0,
\frac{\pi_4(\lambda,k)}{2}\Bigr],\,\,\,\,k^4<\lambda,
$$
where
\begin{equation}
\pi_4(\lambda,k)=2\int_0^{\Bigl(\frac{3}{\lambda-k^4}\Bigr)^{1/4}}\frac{ds}{\Bigl[-\frac{k^2}{3}s^2+\sqrt{\frac{k^4}{9}s^4+
1-\frac{\lambda-k^4}{3}s^4}\,\,\,\Bigr]^{1/2}}.
\end{equation}

Since $F'_{\lambda,k}(x)>0$, then there exists an inverse function
$G_{\lambda,k}: \Bigl[0,
\frac{\pi_4(\lambda,k)}{2}\Bigr]\rightarrow \Bigl[0,
\Bigl(\frac{3}{\lambda-k^4}\Bigr)^{1/4}\Bigr]$  defined by
\begin{equation}
\int_0^{G_{\lambda,k}(x)}\frac{ds}{\Bigl[-\frac{k^2}{3}s^2+\sqrt{\frac{k^4}{9}s^4+
1-\frac{\lambda-k^4}{3}s^4}\,\,\,\Bigr]^{1/2}}=x.
\end{equation}
Finally $G_{\lambda,k}(x)$ is extended to $\mathbb{R}$ in the same
way as  applied to $G(\lambda,x)$ above. By this extension we have
$G_{\lambda,k}(x)= 0 \Leftrightarrow x= n\pi_4(\lambda,k) $. Now,
it follows from (3.8) and (3.11) that $ u(x)= c G_{\lambda,k}(x) $
is the general solution for the following equation
\begin{equation}
k^2(k^2u^2+u'^2)u- \bigl((k^2u^2+u'^2)u'\bigr)'=\lambda u^3,\quad
u(0)=0.
\end{equation}
Thus, the main question is: for what values of $k$ and $\lambda$
are there non-trivial solutions, among $ u(x)= cG_{\lambda,k}(x)
$,
 satisfying the condition $u(1)=0$? The answer to this question
 has almost been given in the above discussion. We summarize these in the following theorem.
\begin{theorem} Let $p=4$ and $k\not=0$. In this case we have:

a) All eigen-parameters $(\lambda,k)$ for problem (3.6) lie in the
parabola $k^4<\lambda$;

b) For every $k\in \mathbb{R}$, the set of all eigenvalues
$\lambda$ of problem (3.6) consists of a sequence of positive
numbers $\lambda_n(k)$ such that $\lambda_n(k)\rightarrow\infty$;

c) There is a number $\lambda^*>0$ such that for each $\lambda$
satisfying $\lambda^*\leq\lambda $, the set of all
eigen-parameters $k$ $(k\geq 0)$ of problem (3.6) consists of a
finite number of eigen-parameters $k_1(\lambda),
k_2(\lambda),...,k_{n(\lambda)}(\lambda)$ which belong to the
interval $[0, \lambda^{1/4}]$. Moreover, $n(\lambda)\rightarrow
\infty$ as $\lambda\rightarrow\infty$;

d) In the case of $\lambda <\lambda^* $  problem (3.6) has only
trivial solution.
\end{theorem}
Proof. a) Let $\lambda\leq k^4$. Then we have
$-\frac{k^2}{3}s^2+\sqrt{\frac{k^4}{9}s^4+
1-\frac{\lambda-k^4}{3}s^4}>0$ for all $s\in \mathbb{R}$. Hence,
it follows from (3.10) that the function $F_{\lambda,k}(x)$ is
well defined, positive function on $[0,+\infty)$ and by definition
so is the function $G_{\lambda,k}(x)$. Since the solution of
(3.11) is given by $ u(x)= cG_{\lambda,k}(x) $, then $u(1)=0$
implies $c=0$. Thus, all eigen-parameters $(\lambda,k)$ for
problem (3.6) lie in the parabola $k^4<\lambda$ (see \cite{hasa10}
for more facts on the localization of the eigen-parameters).

b) Actually, this part has already been proved in Theorem 2.1. In
our case we just present an alternative method for the proof. By
(3.9), for a fixed $k$ the function $\pi_4(\lambda,k)$ satisfies
the following conditions: $\,\lim_{\lambda\rightarrow
k^4+}\pi_4(\lambda,k)=+\infty$, $\lim_{\lambda\rightarrow
+\infty}\pi_4(\lambda,k)=0$ and $\pi_4(\lambda,k)$ is a decreasing
function on $(k^4,+\infty)$. It follows from these properties
that, for each $n\in \mathbb{N}$ the equation $\pi_4(\lambda,k)=
\frac{1}{n}$ has a unique solution $\lambda_n$ and
$\lambda_n\rightarrow +\infty$ as $n\rightarrow \infty$.

c) Now, let us fix $\lambda>0$. Then $\pi_4(\lambda,.)$ is defined
on the interval $[0, \lambda^{1/4}]$. Moreover,
$$
\pi_4(\lambda,0)=
2\int_0^{(\frac{3}{\lambda})^{1/4}}\frac{ds}{\bigl(1-\frac{\lambda}{3}s^4\bigr)^{1/4}},
$$
and $ \lim_{k\rightarrow \lambda^{1/4}-}\pi_4(\lambda,k)=
+\infty.$ Let  $\lambda^*$ be the solution of the equation
$\pi_4(\lambda,0)=1$.  $\pi_p(\lambda,0)$ is a decreasing function
and $\pi_p(\lambda,0)\rightarrow 0$ as $\lambda\rightarrow
+\infty$. Now, clearly the equation $\pi_4(\lambda,k)=
\frac{1}{n}$ has a solution if and only if $\lambda^*\leq\lambda $
and $\frac{1}{n}\in [\pi_4(\lambda), 1]$. These facts prove c) and
d).

\subsection{ On the analytical solutions for arbitrary $p>1$}

In the previous subsection  we have  constructed analytical
solutions to (3.1) for some special cases. Unfortunately, this
method can not be applied to problem (3.1) for arbitrary $p>1$.
Indeed if we try to repeat the same steps from the case $p=4$ then
by using the substitution $v=v(u)$ we can reduce the order of the
equation
$$
k^2|\nabla_k u|^{p-2}u- \bigl(|\nabla_k u|^{p-2}u'\bigr)'= \lambda
|u|^{p-2}u
$$
and by solving the reduced equation as a homogeneous, $1^{st}$
order ordinary differential equation we get
$$
u= \frac{(k^2+ \frac{u'^2}{u^2})^{1/p}\Bigl(-k^2(k^2+
\frac{u'^2}{u^2})^{p/2}-\frac{u'^2(k^2+
\frac{u'^2}{u^2})^{p/2}}{u^2}+ \frac{pu'^2(k^2+
\frac{u'^2}{u^2})^{p/2}}{u^2}+\lambda(k^2+\frac{u'^2}{u^2})\Bigr)^{-\frac{1}{p}}}{c}.
$$
However, this equation is not a radically solvable algebraic
equation with respect to $u'$. Therefore, in this case we have to
use a different method. We study these questions in the next
section.

\section{On the dispersion relation $k(\lambda)$ for ar\-bit\-rary $p>1$: A variational approach}

This Section is devoted to the dispersion relation  $k(\lambda)$,
when the problem is not analytically solvable. In this case we do
not need to restrict ourselves to one dimensional problems. Our
problem is
\begin{equation}
 k^2\int_\Omega|\nabla_k u|^{p-2}uv\,dx+ \int_\Omega|\nabla_k
u|^{p-2}\nabla u.\nabla v\,dx= \lambda\int_\Omega
|u|^{p-2}uv\,dx,\,\,p>1,
\end{equation}
where $\Omega$ is a bounded domain in $\mathbb{R}^n$, $u\in
W_0^{1,p}(\Omega)$, $\nabla_ku:= (ku, \nabla u)$ and $\nabla
u=(u_{x_1}, u_{x_2},...,u_{x_n})$. The approach we use below is
based on the methods of operator pencils. Now we briefly describe
the method of operator pencils which has widely been used in the
spectral theory of linear operator pencils (see \cite{mark88} and
\cite{cola06}). An operator pencil is an operator-valued function
and particularly, it is a polynomial with coefficients in a space
of linear operators.  Typical eigenvalue problems for operator
pencils are:
$$
k^nA_n u+ k^{n-1}A_{n-1}u+...+ kA_1u+A_0u=0
$$
or

$$
k^nA_nu+ k^{n-1}A_{n-1}u+...+ kA_1u+A_0u= \lambda Du.
$$
In the case $k=2$ we deal with a quadratic eigenvalue problem and
we are going to present variational techniques in this case. Thus
we have a quadratic eigenvalue problem

\begin{equation}
 k^2A u+ kBu+Cu= \lambda Du.
\end{equation}
Let us consider the equation
$$
k^2(Au,u)+k(Bu,u)+(Cu,u)= \lambda (Du,u).
$$
This equation defines the following functionals
$$
r_\pm(u,\lambda)= \frac{-(Bu,u) {\pm}\sqrt{(Bu,u)^2 -4((C-\lambda
D)u,u)(Au,u)}}{2(Au,u)}
$$
which play a central role in the variational theory of the
eigenvalue problems (4.2). Actually, in this theory the
functionals $r_\pm(u,\lambda)$ are play the same role as the
functionals $\Phi(u)$ and $G_k(u)$ (see Section 2) in the
Ljusternic-Schnirelman critical point theory. It turns out that
all variational characterizations for $k$ for problem (4.2) are
obtained via $r_\pm(u,\lambda)$. Namely, under some additional
conditions we have (see \cite{cola06} and references therein):

 \begin{equation}
k^\pm_ n(\lambda)=\inf_{L\subset X\atop \dim L=n}\sup_{u\in L\atop
x\neq 0}r_\pm(u,\lambda).
\end{equation}

Now we follow this method for nonlinear problems. We have seen a
simple connection between  problems (4.1) and (4.2) in the case
$p=4$. We have established a connection between $k$ and $\lambda$
by the following differential equations (see the equation next to
(3.6)):
$$
k^4u^3-(uu'^2+u^2u'')k^2-u'^2u''=\lambda u^3,
$$
which can be written in the operator pencil form (we set $k^2:=
\nu$)
\begin{equation}
\nu^2Au+\nu Bu+ Cu=\lambda Au,
\end{equation}
where $Au= u^3$, $Bu=-(uu'^2+u^2u'')$ and $Cu= -u'^2u''$.

 Although the operators $A,\,B$ and $C$ in (4.4) are non-linear  we
 can extend many methods applied in the spectral theory of the operator pencils to the
 non-linear eigenvalue problems,
including problem (4.1) and its  particular case (4.4). Below we
give some results in this direction.

If we replace $k^2$ by $\nu$ in (4.1) then
\begin{equation}
\begin{split}
&\nu\int_\Omega(\nu u^2+|\nabla u|^2)^\frac{p-2}{2}uv\,dx+
\int_\Omega (\nu u^2+|\nabla u|^2)^\frac{p-2}{2}\nabla
u\cdot\nabla
v\,dx\\
&= \lambda\int_\Omega |u|^{p-2}uv\,dx,
\end{split}
\end{equation}
where $\nu\geq 0$. Equation (4.5) is the variational equation for
the functional $F_{\lambda, \nu}(u)= \int_\Omega(\nu u^2+|\nabla
u|^2)^{\frac{p}{2}}\,dx-\lambda\int_\Omega|u|^p\,dx$. Let us fix
$\lambda>0$ and define
 $f_\lambda(\nu,u):= F_{\lambda,\nu}(u)$. Evidently, $f_\lambda:
\mathbb{R}\times X\rightarrow \mathbb{R}$ is  a continuously
differentiable functional. The equation $f_\lambda(\nu,u)=0$ (for
a fixed $\lambda$) defines the so-called root functional
$r_\lambda(u):= \nu$  with values in an interval $J\subset [0,
\sqrt[p]{\lambda^2}]$, which plays the same role as $r_\pm(u,
\lambda)$ for (4.2). We note that the equation
$f_\lambda(\nu,u)=0$ in general defines several functionals and
each functional describes the eigenvalues which belong to its
range.

Next, we give the basic relation between the problem (4.5) and the
root functional $r_\lambda(u)$ (see also \cite{hasa10}). In what
follows we fix $\lambda$ and by eigenvalues we mean a parameter
$\nu$ (or the same $k^2$), satisfying (4.5) with a non-trivial
$u\in W_0^{1,p}(\Omega)$.

\begin{theorem} a) $r(u)\in C^1(X\setminus\{0\},J)$ and it is
extended as a continuous mapping on $X$ by setting $r(0)=0$,

b) $(\lambda, \nu), \,\,\nu\in J$ is an eigen-pair, corresponding
to the eigenvector $u$ for  problem (4.5)  if and only if $u$ is a
critical point and $\nu$ is a critical level for $r$, i.e.,
$\langle r'(u),v\rangle= 0$ for all $v\in X$ and $r(u)=\nu,$

c) All eigenvalues lie in the parabola
$\nu^{\frac{p}{2}}<\lambda$,

d) If $J$ is closed, then the end points of the interval $J$ are
eigenvalues of problem (4.5).

\end{theorem}
Proof.  a) As $$\frac{d}{d\nu}f_\lambda(\nu,u)\Bigl |_{\nu= r(u)}=
\frac{p}{2}\int_\Omega(r(u) u^2+|\nabla
u|^2)^{\frac{p-2}{2}}u^2\,dx>0, \quad 0\not= u\in X ,$$ it follows
from the implicit function theorem that $r\in
C^1(X\setminus\{0\},J)\,$ (see \cite{ze85}, vol I, p.149).

b) By the definition of $r(u)$, we have
$$
\int_\Omega(r(u)u^2+|\nabla u|^2)^{\frac{p}{2}}\,dx=
\lambda\int_\Omega |u|^p\,dx.
$$
Taking the  Fr\'{e}chet derivative from both  sides, we obtain
\begin{equation*}
\begin{split}
&\frac{p}{2}\int_\Omega(r(u)u^2+|\nabla
 u|^2)^{\frac{p-2}{2}}\Bigl[\langle r'(u),v\rangle u^2+ 2
 r(u)uv+2\nabla u\cdot\nabla v\Bigr]\,dx\\
& =\lambda\int_\Omega |u|^{p-2}uv\,dx.
\end{split}
\end{equation*}

%$\frac{p}{2}\int_\Omega(r(u)u^2+|\nabla
% u|^2)^{\frac{p-2}{2}}\Bigl[\langle r'(u),v\rangle u^2+ 2
% r(u)uv+2\nabla u\cdot\nabla v\Bigr]\,dx= \\ \lambda\int_\Omega |u|^{p-2}uv\,dx$.

By regrouping the terms, we get
\begin{equation*}
\begin{split}
&\frac{p}{2}\langle r'(u),v\rangle
\int_\Omega|\nabla_{r(u)}u|^{p-2}u^2\,dx+
p\Bigl[r(u)\int_\Omega|\nabla_{r(u)}u|^{p-2}uv\,dx+\\
&\int_\Omega |\nabla_{r(u)}u|^{p-2}\nabla u\cdot\nabla v\,dx-
\lambda\int_\Omega |u|^{p-2}uv\,dx\Bigr]=0,
\end{split}
\end{equation*}

%$ \frac{p}{2}\langle r'(u),v\rangle
%\int_\Omega|\nabla_{r(u)}u|^{p-2}u^2\,dx+
%p\Bigl[r(u)\int_\Omega|\nabla_{r(u)}u|^{p-2}uv\,dx+\\\int_\Omega
%|\nabla_{r(u)}u|^{p-2}\nabla u\cdot\nabla v\,dx- \lambda\int_\Omega
%|u|^{p-2}uv\,dx\Bigr]=0 $
or
$$
\frac{1}{2}\langle r'(u),v\rangle
\int_\Omega|\nabla_{r(u)}u|^{p-2}u^2\,dx + \langle
F_{\lambda,r(u)}'(u),v\rangle=0.
$$
Finally,
\begin{equation}
\langle r'(u),v\rangle= \frac{-2\langle
F_{\lambda,r(u)}'(u),v\rangle}{\int_\Omega|\nabla_{r(u)}u|^{p-2}u^2\,dx
},\,\,\,\,u\not= 0.
\end{equation}
and by the definition of $F_{\lambda,\nu}'(u)$ a pair of numbers
$(\lambda,\nu)$ is an eigen-pair  if and only if $\langle
F_{\lambda,\nu}'(u),v\rangle= 0$. The needed results follow from
(4.6).

c) This fact immediately follows from the inequality
$$
\nu^{\frac{p}{2}}\int_\Omega |u|^p\,dx< \int_\Omega(\nu
u^2+|\nabla u|^2)^{\frac{p}{2}}=\lambda\int_\Omega|u|^p\,dx,
$$
where $u\not=0$.

d) The closeness of $J$   means that $\inf r(u)$ and $\sup r(u)$
attains. Consequently, these points are critical levels for the
functional $r(u)$.
  $\,\,\,\Box$

Finally, there are a finite number of eigenvalues  for problem
(4.1) denoted by $k_1(\lambda),
k_2(\lambda),...,k_{n(\lambda)}(\lambda)$ (see Theorem 3.1), which
are described by
$$
k_n(\lambda)=\inf_{K\subset \mathcal{K}_n}\sup_{u\in
K}r_\lambda(u).
$$

For this it is enough to check the Palais-Smale condition for
$r_\lambda(u)$ in the interval $J\subset [0, \sqrt[p]{\lambda}]$
or  conditions $H1-H4$ given in Section 2.

 \hfil\break

%\emph{Address:}
% Department of Mathematics, Dogus University, Acibadem, Kadik\"{o}y, 34722, Istanbul,
% Turkey
%
%\emph{e-mail:} hasanov61@yahoo.com;  \hfil  mhasansoy@dogus.edu.tr
\end{document}